# Representation of Solutions of Linear Homogeneous Caputo Fractional Differential Equations with Continuous Variable Coefficients

Sun-Ae PAK , Myong-Ha KIM and Hyong-Chol O *

*Faculty of Mathematics, **Kim Il Sung** University, Pyongyang, D. P. R. Korea*

* Corresponding author email: ohyongchol@yahoo.com

## Abstract

We consider the canonical fundamental systems of solutions of linear homogeneous Caputo fractional differential equations with continuous variable coefficients. Here we gained a series-representation of the canonical fundamental system by coefficients of the considered equations and the representation of solution to initial value problems using the canonical fundamental system. According to our results, the canonical fundamental system of solutions to linear homogeneous differential equation with Caputo fractional derivatives and continuous variable coefficients has different representations according to the distributions of the lowest order of the fractional derivatives in the equation and the distance from the highest order to its adjacent order of the fractional derivatives in the equation.

## 1. Introduction

It is well known that fractional derivative is excellent tool for description of many phenomena (including memory and heredity effects) in fields of science and technology [10]. So recently, many authors are widely studying the theory of fractional differential equations and their application to other fields of science and technology and the topics includes not only the existence, uniqueness and representation of solutions to linear equations but also analytical and numerical methods, stability, bifurcation, chaos and control problems in (systems of) nonlinear fractional equations. (See [11] and another reviewing articles and papers of that theme section,)

In this paper, we study a series representation that analytically represents the canonical fundamental system of solutions of a linear homogeneous Caputo fractional differential equation with continuous variable coefficients only using its coefficients.

One advantage of using fractal calculus in studying differential equations is that it can provide closed formulae (series representation) of solutions even in the case with variable coefficients. In Bonilla et al [1], they obtained an explicit representation of the general solution to the system of linear homogeneous Caputo fractional differential equation with constant coefficients in the case that it includes one term of fractional derivative. In Hu et al [2], they considered linear Caputo fractional differential equation with $n$ terms of fractional derivatives and constant coefficients and obtain solutions of this kind of fractional differential equations by Adomian decomposition method. In Luchko et al [5], they developed the operational calculus of Mikusinŝki's type for Caputo fractional differential operators and applied it to obtain the exact solutions of initial value problems for linear fractional differential equations with constant coefficients and fractional derivatives in Caputo's sense. In Morita et al [6], they obtained solutions of initial value problems of fractional differential equations with constant coefficients using Neumann series for the corresponding Volterra integral equations and expressed by generalized Mittag-Leffler functions.

In [12] an operational calculus of the Mikusinskis type was introduced for the generalized Riemann-Liouville fractional derivative operator and it was applied to solve the corresponding initial value problem for a general $n$-term linear fractional differential equation with constant coefficients and GRLFD of arbitrary orders and types. But they provided initial values independent on the type of generalized Riemann-Liouville fractional derivatives in n-terms linear fractional differential equation and defined the space of solutions independent on the type. They did not provide a necessary and sufficient condition for existence of the solution of n-terms linear fractional differential equation and thus they made some mistakes as mentioned in [9]. In [9], authors present the existence and representation of solution for an $n$-terms linear initial value problem with generalized Riemann-Liouville fractional derivatives and constant coefficients by using operational calculus.

The study on the fractional differential equations with variable coefficients are also provided. In Kilbas et al [3], they investigated solutions around an ordinary point for linear homogeneous Caputo fractional differential equations with sequential fractional derivatives of order $k\alpha(0<\alpha\le 1)$ and variable coefficients. In [8], authors studied on explicit representations of Green's function for linear (Riemann-Liouvilles) fractional differential operators with variable coefficients continuous in $[0,\infty)$ and applied it to obtain explicit representations for solution of non-homogeneous fractional differential equation with variable coefficients of general type.

In [7], authors studied the existence and uniqueness of the solution to a nonlinear differential equation with Caputo fractional derivative in the space of continuously differentiable functions using Banach fixed point theorem. In [13], they studied an approximation method (spline collocation method) to solve nonlinear fractional differential equations with initial conditions or boundary conditions.

In this article we consider the canonical fundamental systems of solutions of a linear homogeneous Caputo fractional differential equation with $(m+1)$ terms and variable coefficients continuous on the interval $[0, T]$. Here we gained a series-representation of the canonical fundamental system by coefficients of the considered equations and the representation of solution to initial value problems using the canonical fundamental system. Our *main method* is the successive approximation. Our *main result* is the discovery of 5 *patterns of distributions of fractional orders* in the equation which *determine the representations* of the solutions.

According to our results, the canonical fundamental system of solutions to linear homogeneous differential equation with Caputo fractional derivatives and continuous variable coefficients has different representations according to the distributions of the lowest order of the fractional derivatives in the equation and the distance from the highest order to its adjacent order of the fractional derivatives in the equation.

The remainder the paper is organized as follows. In section 2 we set our problem and provide some preliminaries including basic concepts and lemmas. Furthermore, we give 3 main patterns of distributions of fractional orders of the considered equation which determine the representations of the solutions. In sections 3, 4 and 5 we provide series-representations of the canonical fundamental systems by coefficients of the considered equations and the representation of solution to initial value problems using the canonical fundamental system in those 3 cases, respectively. In the last section we give some conclusions.

## 2. Problems and Lemmas (Preliminaries)

In this paper we study the following linear homogeneous Caputo fractional differential equation with variable coefficients:

$$ {}^{c}D_{0+}^{\alpha_0} y(t) = -\sum_{i=1}^{m} a_i(t) \cdot {}^{c}D_{0+}^{\alpha_i} y(t), \quad 0 < t < T \ . \tag{1} $$

Here $T > 0$ is an arbitrary positive number, $\alpha_i (i = 0, 1, \cdots, m)$ are nonnegative real numbers such that $\alpha_0 > \alpha_1 > \cdots > \alpha_m \geq 0$ and there exist nonnegative integers $n_i (i = 0, 1, \cdots, m)$ such that $n_i - 1 < \alpha_i \leq n_i \ (i = 0, 1, \cdots, m)$. It is evident that $n_0 \geq n_1 \geq \cdots \geq n_m \geq 0$.

***Remark 1***. The integer $n_0$ dependent on the highest order $\alpha_0$ is important and it provides the number of initial conditions to give in order to ensure the existence and uniqueness of the solution to the initial value problem of (1).

Now we study an initial value problem of the equation (1) with the following initial conditions:

$$ D^k y(t)\Big|_{t=+0} = b_k \in R, \ k = 0, 1, \cdots, n_0 - 1 . \tag{2} $$

We use the function spaces $C_\gamma[a, b]$ and $C_\gamma^n[a, b]$ introduced in [2]. When $0 \leq \gamma < 1$, we denote by $C_\gamma[a,b]$ the space of all functions $f$ that $t \mapsto (t-a)^\gamma f(t) \in C[a, b]$ and $f$ itself is defined on $(a,b]$. In particular $C_0[a,b] = C[a,b]$. For $n \in N$ we denote by $C_\gamma^n[a,b]$ the space of functions $f$ which are $n-1$ times continuously differentiable on $[a,b]$ and have the derivative $f^{(n)}(x)$ such that $f^{(n)}(x) \in C_\gamma[a,b]$. And we denote as follows:

$$ C_\gamma^0[a,b] = C_\gamma[a,b]. $$

(The definition of their norms of these spaces are also provided in [7].)

**Definition 1**[2] Let $\alpha \in R_+$, $f \in C_\gamma[0, T]$, $0 \le \gamma < 1$. The *fractional integrals* $I_{0+}^{\alpha} f$ of order $\alpha\,(\alpha > 0)$ of function $f$ in the means of Riemann-Liouville are defined by

$$I_{0+}^{\alpha} f(t) = \frac{1}{\Gamma(\alpha)} \int_0^t (t-\tau)^{\alpha-1} f(\tau) d\tau \,, t > 0 .$$

In particular we denote as $I^0 f(t) = f(t)$.

**Definition 2**[2] When $n \in N$, $n-1 < \alpha \le n$, the *Caputo fractional derivatives* ${}^cD_{0+}^{\alpha} f$ of order $\alpha\,(\alpha > 0)$ of function $f$ are defined by

$${}^cD_{0+}^{\alpha} f(t) = D_{0+}^{\alpha}\left[ f(t) - \sum_{k=0}^{n-1} \frac{f^{(k)}(0)}{k!} t^k \right].$$

Here $D_{0+}^{\alpha}$ is *Riemann-Liouville fractional derivative*:

$$D_{0+}^{\alpha} f(t) := \frac{d^n}{dt^n} I_{0+}^{n-\alpha} f(t) \;,\; t > 0 .$$

**Lemma 1.**[4] Let $n \in \mathbf{N}_0 = \{0, 1, \cdots\}$ and $\gamma \in \mathbf{R}$ $(0 \le \gamma \le 1)$. The space $C_\gamma^n[a,b]$ consists of those and only those functions $f$ which are represented in the form

$$f(t) = \frac{1}{(n-1)!} \int_a^t (t-\tau)^{n-1} \varphi(\tau) d\tau + \sum_{k=1}^{n-1} C_k (t-a)^k \;,$$

where $\varphi \in C_\gamma[a, b]$ and $C_k$ $(k = 0, 1, \cdots, n-1)$ are arbitrary constants.

**Definition 3.** A system of functions $y_j(t)(j = 0, 1, \cdots, n_0 - 1)$ is called a ***canonical fundamental system*** of solutions of the homogeneous equation (1) if it satisfies

$${}^cD_{0+}^{\alpha_0} y_j(t) = -\sum_{i=1}^{m} a_i(t) {}^cD_{0+}^{\alpha_i} y_j(t), \quad 0 < t < T \;,$$

$$D^k y_j(t)\Big|_{t=+0} = \begin{cases} 1, & j = k, \\ 0, & j \ne k;\; k, j = 0, 1, \cdots, n_0 - 1. \end{cases}$$

The ***main conclusion*** of this paper is that the solutions to (1) have the *different representations* according to the *different distributions* of fractional orders $\alpha_i\,(i = 0, 1, \cdots, m)$. In order to explain this situation clearly, we introduce the following index sets $H_j$ $(i = 0, 1, \cdots, m)$ to help to show the distribution of given orders $\alpha_i\,(i = 0, 1, \cdots, m)$.

$$H_j = \{i : 0 \le \alpha_i \le j, i = 1, \cdots, m\}, \;\; j = 0, 1, \cdots, n_0 - 1, \tag{3}$$

***Remark 2***. $k \in H_j \Rightarrow \alpha_k \le j$ and $H_i \subset H_j (i < j)$. If $H_j \ne \phi$ (where $\phi$ is the empty set), we let $h_j = \min H_j$ which is the smallest index of fractional orders of (1) that do not exceed $j$. Then it is evident that $m - h_j + 1$ is the number of elements of $H_j$ and $h_i \ge h_j (i < j)$. Thus when $i < j$, $h_i - h_j$ represents the number of such orders $\alpha_k$ that $i < \alpha_k \le j$ and in particular, if $h_i = h_j$ then (1) has no such fractional orders $\alpha_k$ that $i < \alpha_k \le j$.

We will consider the following 3 possible cases which represent different distributions of fractional orders. This is the *main idea* of this paper which allows us successfully give representations of solutions to (1).

**Case I**: $H_0 \neq \phi$. In this case $H_j \neq \phi, \forall j = 0, 1, \cdots, n_0 - 1$ and our equation (1) has the linear term of unknown function $y(t)$ (the term of 0-th derivative). Then (1) is called the ***type I***. In this case we have $\alpha_m = 0,\ h_0 = m,\ n_m = 0$ and

$$n_m = \alpha_m = \alpha_{h_0} = 0 < \alpha_{m-1} < \cdots < \alpha_{h_1} \leq 1 < \alpha_{h_1 - 1} < \cdots < \alpha_{h_j} \leq j < \alpha_{h_j - 1} < \cdots < \alpha_{h_{n_0}} = \alpha_0 \leq n_0. \tag{*}$$

**Case II**: $n_0 \geq 2$ and there exists a $j_0 \in \{0, 1, \cdots, n_0 - 2\}$ such that $H_{j_0} = \phi$ and $H_{j_0 + 1} \neq \phi$. In this case $H_j = \phi,\ j = 0, 1, \cdots, j_0$ and $H_j \neq \phi,\ j = j_0 + 1, \cdots, n_0 - 1$. That is, there is no the term of derivative of order of $j_0$ or lower than $j_0$ in the equation (1). Then (1) is called the ***type II***. In this case we have $\alpha_m \in (0,\ n_0 - 1]$, $1 \leq n_m \leq n_0 - 1$, $j_0 = n_m - 1$ and

$$0 \leq n_m - 1 = j_0 < \alpha_m < \alpha_{m-1} < \cdots < \alpha_{h_{j_0 + k}} \leq j_0 + k < \alpha_{h_{j_0 + k} - 1} < \cdots < \alpha_{h_{n_0 - 1}} \leq n_0 - 1 < \alpha_{h_{n_0 - 1} - 1} < \cdots < \alpha_{h_{n_0}} = \alpha_0 \leq n_0. \tag{**}$$

**Case III**: $H_{n_0 - 1} = \phi$. In this case $H_j = \phi,\ j = 0,\ 1, \cdots,\ n_0 - 1$. That is, there only are the terms of derivative of order of higher than $n_0 - 1$ in the equation (1). Then (1) is called the ***type III***. In this case $\alpha_m \in (n_0 - 1,\ n_0]$ and thus $\alpha_i \in (n_0 - 1,\ n_0]$ for all $i = 0, 1, \cdots, m$. From the definition 2, we know that in this case our equation (1) becomes a very simple integro-differential equation which can be transformed in to an integral equation with nonhomogeneous term of a polynomial on *t* with the degree of $n_0 - 1$.

Now let consider the following approximation sequence for (1), (2):

$$y^0(t) = \sum_{k=0}^{n_0 - 1} b_k \Phi_{k+1}(t)\ , \quad \Phi_{k+1}(t) = t^k / k!,$$

$$y^l(t) = y^0(t) - I_{0+}^{\alpha_0}\left[\sum_{i=1}^{m} a_i(t) \cdot {}^c D_{0+}^{\alpha_i} y^{l-1}(t)\right],\ l = 1,\ 2, \cdots. \tag{4}$$

Then we have the following lemma.

**Lemma 2.** (i) *Assume that* $n_0 > n_1$ *and* $a_i \in C[0, T]$, $i = 1, \cdots, m$. *Then The initial value problem (1), (2) has a unique solution* $y(t) \in C^{\alpha_0, n_0 - 1}[0, T]$, *which is the limit in* $C^{\alpha_0, n_0 - 1}[0, T]$ *of the approximation sequence (4).*

(ii) *Let* $\gamma = n_0 - \alpha_0 (\in [0, 1))$. *Assume that* $n_0 = n_1$ *and. If* $a_i \in C_\gamma^1[0, T]$, $D_{0+}^\gamma a_i \in C[0, T]$, $i = 1, \cdots, m$, *then the initial value problem (1), (2) has a unique solution* $y(t) \in C_\gamma^{n_0}[0, T]$, *which is the limit in* $C_\gamma^{n_0}[0, T]$ *of the approximation sequence (4).*

**Proof:** Let $f(t, {}^cD_{0+}^{\alpha_1}y(t),\cdots, {}^cD_{0+}^{\alpha_m}y(t)) = -\sum_{i=1}^{m} a_i(t) {}^cD_{0+}^{\alpha_i}y(t)$. Then by assumption we have $a_i \in C[0,T]$, $i=1,\cdots,m$ or $0 \le \exists\gamma < 1$ : $a_i \in C_\gamma^1[0,T], \mathrm{i}=1,\cdots,m$ and thus $f:[0,T]\times Q \to R$ ($Q \subset R^m$) has continuity and Lipchitz condition. By the result of [7], the lemma 2 is proved. (QED)

***Remark 3.*** If $0 \le \exists\gamma < 1 : a_i \in C_\gamma^1[0,T],\ i=1,\cdots,m$, then the solution $y(t)$ of the initial value problem (1), (2) satisfies $y(t) \in C^{n_0-1}[0,T] \cap C^{n_0}(0,T]$.

## 3. The Canonical Fundamental System of Solutions to the Equations of the Type I

**Theorem 1.** *Let $\gamma = n_0 - \alpha_0 (\in [0,1))$ and Assume that $n_0 = n_1$ and $a_i \in C_\gamma^1[0,T]$, $D_{0+}^{n_0-\alpha_0} a_i \in C[0,T], i=1,\cdots,m$. Assume that $H_0 \ne \phi$ (that is, (1) is the type I). Then there exists the unique canonical fundamental system $y_j(t) \in C_\gamma^{n_0}[0,T]$, ($j=0,1,\cdots,\ n_0-1$ of solutions to (1) and it is written by*

$$y_j(t) = \Phi_{j+1}(t) + \sum_{k=0}^{\infty} (-1)^{k+1} I_{0+}^{\alpha_0} \left[ \sum_{i=1}^{m} a_i(t) I_{0+}^{\alpha_0-\alpha_i} \right]^k \sum_{i=h_j}^{m} a_i(t) \Phi_{j+1-\alpha_i}(t), \tag{5}$$

$$j = 0,1,\cdots,n_0-1.$$

**Proof.** From lemma 2-(ii) the existence and uniqueness of canonical fundamental systems is evident. Now let find the canonical system as the limit in $C_\gamma^{n_0}[0,T]$ of the approximation sequence

$$y_j^0(t) = \Phi_{j+1}(t), \tag{6}$$

$$y_j^{l+1}(t) = \Phi_{j+1}(t) - I_{0+}^{\alpha_0} \left[ \sum_{i=1}^{m} a_i(t) {}^cD_{0+}^{\alpha_i} y_j^l(t) \right], \quad l = 0,1,2,\cdots \tag{7}$$

First, we find $y_0(t)$. Let $j=0$ in (6), then $y_0^0(t) = \Phi_1(t)$ and from (7) we have

$$y_0^1(t) = \Phi_1(t) - I_{0+}^{\alpha_0} \sum_{i=1}^{m} a_i(t) {}^cD_{0+}^{\alpha_i} y_0^0(t) = \Phi_1(t) - I_{0+}^{\alpha_0} \sum_{i=1}^{m} a_i(t) {}^cD_{0+}^{\alpha_i} \Phi_1(t). \tag{8}$$

Here we calculate ${}^cD_{0+}^{\alpha_i}\Phi_1(t)$. By the definition 2 we have

$${}^cD_{0+}^{\alpha_i}\Phi_1(t) = D_{0+}^{\alpha_i} \left[ \Phi_1(t) - \sum_{k=0}^{n_i-1} D^k \Phi_1(0) \Phi_{k+1}(t) \right],\ i = 1,\cdots,m. \tag{9}$$

Since $H_0 \ne \phi$, we have $h_0 = \min H_0 = m$. For $i = m$, $\alpha_m = n_m = 0$ and thus

$${}^cD_{0+}^{\alpha_i}\Phi_1(t) = D_{0+}^{\alpha_m}\Phi_1(t) = \Phi_1(t).$$

If $i=1,\cdots,m-1$ then $\alpha_i>0$ and thus from $n_i-1<\alpha_i\le n_i$ we have $n_i\ge 1$. Since

$$D^k\Phi_1(0)=\begin{cases}1, & k=0,\\ 0, & k\ne 0,\end{cases}$$

in (9), then ${}^cD_{0+}^{\alpha_i}\Phi_1(t)=0$ for $i=1,\cdots,m-1$ and

$${}^cD_{0+}^{\alpha_i}\Phi_1(t)=\begin{cases}\Phi_1(t), & i=h_0=m,\\ 0, & i=1,\cdots,m-1.\end{cases}$$

Substituting this into (8), then the first approximation of $y_0(t)$ is given by

$$y_0^1(t)=\Phi_1(t)-I_{0+}^{\alpha_0}a_m(t)\Phi_1(t) \tag{10}$$

and since $h_0=m$ and $\alpha_m=0$, we can rewrite it as follows:

$$y_0^1(t)=\Phi_1(t)-I_{0+}^{\alpha_0}\sum_{i=h_0}^{m}a_i(t)\Phi_{1-\alpha_i}(t).$$

Thus we get the first term( $k=0$ ) of (5) in the case of $j=0$.

Now prove $y_0^1(t)\in C_\gamma^{n_0}[0,T]$. From (10), we have

$$D^{n_0-1}y_0^1(t)=D^{n_0-1}\Phi_1(t)-D^{n_0-1}I_{0+}^{\alpha_0}a_m(t)\Phi_1(t)=-I_{0+}^{\alpha_0-n_0+1}a_m(t)\Phi_1(t).$$

From the assumption $a_m(t)\Phi_1(t)=a_m(t)\in C[0,T]$ and $\alpha_0-n_0+1>0$, and thus we have $I_{0+}^{\alpha_0-n_0+1}a_m(t)\Phi_1(t)\in C[0,T]$, that is,

$$y_0^1(t)\in C^{n_0-1}[0,T].$$

On the other hand, $a_m(t)\Phi_1(t)=a_m(t)\in C_\gamma^1[0,T]$ and $\alpha_0-n_0+1>0$, and thus by lemma1 we have $I_{0+}^{\alpha_0-n_0+1}a_m(t)\Phi_1(t)\in C_\gamma^1[0,T]$. Therefore

$$D^{n_0}y_0^1(t)=D^{n_0}\Phi_1(t)-D^{n_0}I_{0+}^{\alpha_0}a_m(t)\Phi_1(t)=-DI_{0+}^{\alpha_0-n_0+1}a_m(t)\Phi_1(t)\in C_\gamma[0,T].$$

Thus we proved $y_0^1(t)\in C_\gamma^{n_0}[0,T]$.

Next we consider the case that $j=0,\ l=1$ in (7) to find the second approximation of $y_0(t)$.

$$y_0^2(t)=\Phi_1(t)-I_{0+}^{\alpha_0}\sum_{i=1}^{m}a_i(t)^cD_{0+}^{\alpha_i}y_0^1(t)=$$

$$=\Phi_1(t)-I_{0+}^{\alpha_0}\sum_{i=1}^{m}a_i(t)^cD_{0+}^{\alpha_i}\left[\Phi_1(t)-I_{0+}^{\alpha_0}\sum_{i=h_0}^{m}a_i(t)\Phi_{1-\alpha_i}(t)\right]$$

$$=\Phi_1(t)-I_{0+}^{\alpha_0}\sum_{i=1}^{m}a_i(t)^cD_{0+}^{\alpha_i}\Phi_1(t)+I_{0+}^{\alpha_0}\sum_{i=1}^{m}a_i(t)^cD_{0+}^{\alpha_i}I_{0+}^{\alpha_0}\sum_{i=h_0}^{m}a_i(t)\Phi_{1-\alpha_i}(t).$$

Now let $f(t):=\sum_{i=h_0}^{m}a_i(t)\Phi_{1-\alpha_i}(t)$ and calculate ${}^cD_{0+}^{\alpha_i}I_{0+}^{\alpha_0}f(t)$. From the definition 2

$$ {}^cD_{0+}^{\alpha_i} I_{0+}^{\alpha_0} f(t) = D_{0+}^{\alpha_i}\left[ I_{0+}^{\alpha_0} f(t) - \sum_{k=0}^{n_i-1} D^k I_{0+}^{\alpha_0} f(t) \right]_{t=+0} \Phi_{k+1}(t). $$

Here since $k \le n_i - 1 < \alpha_i < \alpha_0$, $i = 1, \cdots, m$ we have $\alpha_0 - k > 0$, $k = 0, 1, \cdots, n_i - 1$ and thus we can rewrite as follows:

$$ D^k I_{0+}^{\alpha_0} f(t) = D^k I^k I_{0+}^{\alpha_0 - k} f(t) = I_{0+}^{\alpha_0 - k} f(t). $$

Since $f(t) \in C[0, T]$ we have

$$ D^k I_{0+}^{\alpha_0} f(t)\big|_{t=+0} = I_{0+}^{\alpha_0 - k} f(t)\big|_{t=+0} = 0, $$

and therefore

$$ {}^cD_{0+}^{\alpha_i} I_{0+}^{\alpha_0} f(t) = D_{0+}^{\alpha_i} I_{0+}^{\alpha_0} f(t) = D_{0+}^{\alpha_i} I_{0+}^{\alpha_i} I_{0+}^{\alpha_0 - \alpha_i} f(t) = I_{0+}^{\alpha_0 - \alpha_i} f(t). $$

Thus the second approximation is given as follows:

$$ \begin{aligned} y_0^2(t) &= \Phi_1(t) - I_{0+}^{\alpha_0} \sum_{i=1}^{m} a_i(t)\, {}^cD_{0+}^{\alpha_i} \Phi_1(t) + I_{0+}^{\alpha_0} \sum_{i=1}^{m} a_i(t) I_{0+}^{\alpha_0 - \alpha_i} \sum_{i=h_0}^{m} a_i(t) \Phi_{1-\alpha_i}(t) \\ &= \Phi_1(t) - I_{0+}^{\alpha_0} a_m(t) \Phi_1(t) + I_{0+}^{\alpha_0} \sum_{i=1}^{m} a_i(t) I_{0+}^{\alpha_0 - \alpha_i} \sum_{i=h_0}^{m} a_i(t) \Phi_{1-\alpha_i}(t) \\ &= \Phi_1(t) + (-1)^1 I_{0+}^{\alpha_0} \sum_{i=h_0}^{m} a_i(t) \Phi_{1-\alpha_i}(t) + (-1)^2 I_{0+}^{\alpha_0} \sum_{i=1}^{m} a_i(t) I_{0+}^{\alpha_0 - \alpha_i} \sum_{i=h_0}^{m} a_i(t) \Phi_{1-\alpha_i}(t) \\ &= \Phi_1(t) + \sum_{k=0}^{1} (-1)^{k+1} I_{0+}^{\alpha_0} \left[ \sum_{i=1}^{m} a_i(t) I_{0+}^{\alpha_0 - \alpha_i} \right]^k \sum_{i=h_0}^{m} a_i(t) \Phi_{1-\alpha_i}(t). \end{aligned} $$

Since $h_0 = m, \alpha_m = 0$, then $\Phi_{1-\alpha_i}(t) = \Phi_1(t) \in C[0, T]$ and therefore

$$ \begin{aligned} D^{n_0} y_0^2(t) &= D^{n_0} \Phi_1(t) + \sum_{k=0}^{1} (-1)^{k+1} D^{n_0} I_{0+}^{\alpha_0} \left[ \sum_{i=1}^{m} a_i(t) I_{0+}^{\alpha_0 - \alpha_i} \right]^k \sum_{i=h_0}^{m} a_i(t) \Phi_{1-\alpha_i}(t) \\ &= \sum_{k=0}^{1} (-1)^{k+1} D\, I_{0+}^{\alpha_0 - n_0 + 1} \left[ \sum_{i=1}^{m} a_i(t) I_{0+}^{\alpha_0 - \alpha_i} \right]^k \sum_{i=h_0}^{m} a_i(t) \Phi_{1-\alpha_i}(t). \end{aligned} $$

From the assumption and lemma 1, we have

$$ I_{0+}^{\alpha_0 - n_0 + 1} \left[ \sum_{i=1}^{m} a_i(t) I_{0+}^{\alpha_0 - \alpha_i} \right]^k \sum_{i=h_0}^{m} a_i(t) \Phi_{1-\alpha_i}(t) \in C_\gamma^1[0, T]. $$

Thus we have $y_0^2(t) \in C_\gamma^{n_0}[0, T]$.

Now under the assumption that the $l$-th approximation of $y_0(t)$ is provided by

$$ y_0^l(t) = \Phi_1(t) + \sum_{k=0}^{l-1} (-1)^{k+1} I_{0+}^{\alpha_0} \left[ \sum_{i=1}^{m} a_i(t) I_{0+}^{\alpha_0 - \alpha_i} \right]^k \sum_{i=h_0}^{m} a_i(t) \Phi_{1-\alpha_i}(t) $$

and $y_0^l(t) \in C_\gamma^{n_0}[0, T]$, we find the $(l+1)$-st approximation of $y_0(t)$.

$$y_0^{l+1}(t) = \Phi_1(t) - I_{0+}^{\alpha_0} \sum_{i=1}^{m} a_i(t)^c D_{0+}^{\alpha_i} y_0^l(t) =$$

$$= \Phi_1(t) - I_{0+}^{\alpha_0} \sum_{i=1}^{m} a_i(t)^c D_{0+}^{\alpha_i} \left\{ \Phi_1(t) + \sum_{k=0}^{l-1} (-1)^{k+1} I_{0+}^{\alpha_0} \left[ \sum_{i=1}^{m} a_i(t) I_{0+}^{\alpha_0-\alpha_i} \right]^k \sum_{i=h_0}^{m} a_i(t)\Phi_{1-\alpha_i}(t) \right\}$$

$$= \Phi_1(t) - I_{0+}^{\alpha_0} \sum_{i=1}^{m} a_i(t)^c D_{0+}^{\alpha_i} \Phi_1(t) -$$

$$- I_{0+}^{\alpha_0} \sum_{i=1}^{m} a_i(t)^c D_{0+}^{\alpha_i} \sum_{k=0}^{l-1} (-1)^{k+1} I_{0+}^{\alpha_0} \left[ \sum_{i=1}^{m} a_i(t) I_{0+}^{\alpha_0-\alpha_i} \right]^k \sum_{i=h_0}^{m} a_i(t)\Phi_{1-\alpha_i}(t)$$

$$= \Phi_1(t) - I_{0+}^{\alpha_0} \sum_{i=h_0}^{m} a_i(t)\Phi_{1\text{-}\alpha_i}(t)$$

$$- I_{0+}^{\alpha_0} \sum_{i=1}^{m} a_i(t) \sum_{k=0}^{l-1} (-1)^{k+1} I_{0+}^{\alpha_0-\alpha_i} \left[ \sum_{i=1}^{m} a_i(t) I_{0+}^{\alpha_0-\alpha_i} \right]^k \sum_{i=h_0}^{m} a_i(t)\Phi_{1-\alpha_i}(t)$$

$$= \Phi_1(t) + (-1) I_{0+}^{\alpha_0} \sum_{i=h_0}^{m} a_i(t)\Phi_{1\text{-}\alpha_i}(t) +$$

$$+ I_{0+}^{\alpha_0} \sum_{k=0}^{l-1} (-1)^{k+2} \sum_{i=1}^{m} a_i(t) I_{0+}^{\alpha_0-\alpha_i} \left[ \sum_{i=1}^{m} a_i(t) I_{0+}^{\alpha_0-\alpha_i} \right]^k \sum_{i=h_0}^{m} a_i(t)\Phi_{1-\alpha_i}(t)$$

$$= \Phi_1(t) + \sum_{k=0}^{l} (-1)^{k+1} I_{0+}^{\alpha_0} \left[ \sum_{i=1}^{m} a_i(t) I_{0+}^{\alpha_0-\alpha_i} \right]^k \sum_{i=h_0}^{m} a_i(t)\Phi_{1-\alpha_i}(t).$$

Thus the $(l+1)$-st approximation of $y_0(t)$ is provided by

$$\mathrm{y}_{0+}^{l+1}(\mathrm{t}) = \Phi_1(t) + \sum_{k=0}^{l} (-1)^{k+1} I_{0+}^{\alpha_0} \left[ \sum_{i=1}^{m} a_i(t) I_{0+}^{\alpha_0-\alpha_i} \right]^k \sum_{i=h_0}^{m} a_i(t)\Phi_{1-\alpha_i}(t).$$

In the similar way as the above we get $y_0^{l+1}(t) \in C_\gamma^{n_0}[0,T]$.

By induction we proved that for any $n = 0, 1, \cdots$, $n$ th approximation of $y_0(t)$ is provided by

$$\mathrm{y}_0^n(\mathrm{t}) = \Phi_1(t) + \sum_{k=0}^{n-1} (-1)^{k+1} I_{0+}^{\alpha_0} \left[ \sum_{i=1}^{m} a_i(t) I_{0+}^{\alpha_0-\alpha_i} \right]^k \sum_{i=h_0}^{m} a_i(t)\Phi_{1-\alpha_i}(t).$$

and $y_0^n(t) \in C_\gamma^{n_0}[0,T]$. From lemma2-(ii) the sequence $\{ y_0^n(t) \}$ converges in $C_\gamma^{n_0}[0,T]$ and we get the first element $y_0(t)$ of canonical system:

$$y_0(t) = \lim_{n\to\infty} y_{0+}^n(\mathrm{t})$$

$$=\Phi_1(t)+\sum_{k=0}^{\infty}(-1)^{k+1}I_{0+}^{\alpha_0}\left[\sum_{i=1}^{m}a_i(t)I_{0+}^{\alpha_0-\alpha_i}\right]^k\sum_{i=h_0}^{m}a_i(t)\Phi_{1-\alpha_i}(t)\in C_\gamma^{n_0}[0,T].$$

Now for any $j(j=1,\cdots,n_0-1)$, we find the $j$-th element $y_j(t)$ of the canonical system. From (6) $y_j^0(t)=\Phi_{j+1}(t)$ and by (7) the first approximation of $y_j(t)$ is given by

$$y_j^1(t)=\Phi_{j+1}(t)-I_{0+}^{\alpha_0}\sum_{i=1}^{m}a_i(t)^cD_{0+}^{\alpha_i}y_j^0(t)=\Phi_{j+1}(t)-I_{0+}^{\alpha_0}\sum_{i=1}^{m}a_i(t)^cD_{0+}^{\alpha_i}\Phi_{j+1}(t).$$

By the definition 2, ${}^cD_{0+}^{\alpha_i}\Phi_{j+1}(t)$ is as follows

$${}^cD_{0+}^{\alpha_i}\Phi_{j+1}(t)=D_{0+}^{\alpha_i}\left[\Phi_{j+1}(t)-\sum_{k=0}^{n_i-1}D^k\Phi_{j+1}(0)\Phi_{k+1}(t)\right]\ (i=1,\cdots,m)..$$

If $i\ge h_j$ then $0\le\alpha_i\le j$ and $n_i\le j$, and thus $D^k\Phi_{j+1}(0)=0,k=0,\cdots,n_0-1$, and thus we have

$${}^cD_{0+}^{\alpha_i}\Phi_{j+1}(t)=D_{0+}^{\alpha_i}\Phi_{j+1}(t),\ i=h_j,\cdots,m.$$

.

If $i<h_j$, then $n_i>j$, so

$$\sum_{k=0}^{n_i-1}D^k\Phi_{j+1}(0)\Phi_{k+1}(t)=\Phi_{j+1}(t)$$

and thus ${}^cD_{0+}^{\alpha_i}\Phi_{j+1}(t)=0$, $i<h_j$. That is, we have

$${}^cD_{0+}^{\alpha_i}\Phi_{j+1}(t)=\begin{cases}D_{0+}^{\alpha_i}\Phi_{j+1}(t), & h_j\le i\le m,\\ 0, & 1\le i<h_j.\end{cases}$$

So the first approximation of $y_j(t)$ is provided by

$$y_j^1(t)=\Phi_{j+1}(t)-I_{0+}^{\alpha_0}\sum_{i=1}^{m}a_i(t)^cD_{0+}^{\alpha_i}\Phi_{j+1}(t)=\Phi_{j+1}(t)-I_{0+}^{\alpha_0}\sum_{i=h_j}^{m}a_i(t)D_{0+}^{\alpha_i}\Phi_{j+1}(t)$$

$$=\Phi_{j+1}(t)-I_{0+}^{\alpha_0}\sum_{i=h_j}^{m}a_i(t)\Phi_{j+1-\alpha_i}(t)$$

and we have $y_j^1(t)\in C_\gamma^{n_0}[0,T]$. Here

$$\Phi_{j+1-\alpha_i}(t):=D_{0+}^{\alpha_i}\Phi_{j+1}(t).$$

The second approximation of $y_j(t)$ is given by

$$y_j^2(t)=\Phi_{j+1}(t)-I_{0+}^{\alpha_0}\sum_{i=1}^{m}a_i(t)^cD_{0+}^{\alpha_i}y_j^1(t)=$$

$$=\Phi_{j+1}(t)-I_{0+}^{\alpha_0}\sum_{i=1}^{m}a_i(t)^cD_{0+}^{\alpha_i}\left[\Phi_{j+1}(t)-I_{0+}^{\alpha_0}\sum_{i=h_j}^{m}a_i(t)\Phi_{j+1-\alpha_i}(t)\right]$$

$$=\Phi_{j+1}(t)-I_{0+}^{\alpha_0}\sum_{i=1}^{m}a_i(t)^cD_{0+}^{\alpha_i}\Phi_{j+1}(t)+I_{0+}^{\alpha_0}\sum_{i=1}^{m}a_i(t)^cD_{0+}^{\alpha_i}I_{0+}^{\alpha_0}\sum_{i=h_j}^{m}a_i(t)\Phi_{j+1-\alpha_i}(t)\,.$$

Here note that ${}^cD_{0+}^{\alpha_i}I_{0+}^{\alpha_0}\sum_{i=h_j}^{m}a_i(t)\Phi_{j+1-\alpha_i}(t)=I_{0+}^{\alpha_0-\alpha_i}\sum_{i=h_j}^{m}a_i(t)\Phi_{j+1-\alpha_i}(t)$, we have

$$y_j^2(t)=\Phi_{j+1}(t)+\sum_{k=0}^{1}(-1)^{k+1}I_{0+}^{\alpha_0}\left[\sum_{i=1}^{m}a_i(t)I_{0+}^{\alpha_0-\alpha_i}\right]^k\sum_{i=h_j}^{m}a_i(t)\Phi_{j+1-\alpha_i}(t)\,.$$

By induction, the $n$-th approximation of $y_j(t)$ is given by

$$y_j^n(t)=\Phi_{j+1}(t)+\sum_{k=0}^{n-1}(-1)^{k+1}I_{0+}^{\alpha_0}\left[\sum_{i=1}^{m}a_i(t)I_{0+}^{\alpha_0-\alpha_i}\right]^k\sum_{i=h_j}^{m}a_i(t)\Phi_{j+1-\alpha_i}(t)\,.$$

Therefore

$$y_j(t)=\lim_{n\to\infty}y_j^n(t)=\Phi_{j+1}(t)+\sum_{k=0}^{\infty}(-1)^{k+1}I_{0+}^{\alpha_0}\left[\sum_{i=1}^{m}a_i(t)I_{0+}^{\alpha_0-\alpha_i}\right]^k\sum_{i=h_j}^{m}a_i(t)\Phi_{j+1-\alpha_i}(t)\,.$$

(QED)

**Corollary 1.** *Under the assumption of theorem 1, a solution* $y(t)\in C_\gamma^{n_0}[0,T]$ *to the initial value problem (1) and (2) uniquely exists and represents by*

$$y(t)=\sum_{j=0}^{n_0-1}b_jy_j(t)\,.$$

*Here* $y_j(t)\in C_\gamma^{n_0}[0,T],\ j=0,1,\cdots,n_0-1$ *is the canonical fundamental system of (1) given by (5).*

**Theorem 2.** *Assume that* $n_0>n_1$, $a_i\in C[0,T]\ (i=1,\cdots,m)$ *and* $H_0\neq\phi$. *Then there exists the unique canonical fundamental system*

$$y_j(t)\in C^{\alpha_0,n_0-1}[0,T],\ j=0,1,\cdots,n_0-1$$

*of solutions to (1) and it is written by*

$$y_j(t)=\Phi_{j+1}(t)+\sum_{k=0}^{\infty}(-1)^{k+1}I_{0+}^{\alpha_0}\left[\sum_{i=1}^{m}a_i(t)I_{0+}^{\alpha_0-\alpha_i}\right]^k\sum_{i=h_j}^{m}a_i(t)\Phi_{j+1-\alpha_i}(t),\qquad(11)$$

$$j=0,1,\cdots,n_1-1,$$

$$y_j(t)=\Phi_{j+1}(t)+\sum_{k=0}^{\infty}(-1)^{k+1}I_{0+}^{\alpha_0}\left[\sum_{i=1}^{m}a_i(t)I_{0+}^{\alpha_0-\alpha_i}\right]^k\sum_{i=1}^{m}a_i(t)\Phi_{j+1-\alpha_i}(t),\qquad(12)$$

$$j=n_1,n_1+1,\cdots,n_0-1.$$

**Proof**. From lemma 2-(i) the existence and uniqueness of canonical systems is evident. Now let find the canonical system as the limit of the approximation sequence (6) and (7). Fixed $j=0,1,\cdots,n_0-1$, we find the $j$-th element $y_j(t)$ of the canonical system.

$$y_j^0(t)=\Phi_{j+1}(t),$$

$$y_j^1(t)=\Phi_{j+1}(t)-I_{0+}^{\alpha_0}\sum_{i=1}^{m}a_i(t)^cD_{0+}^{\alpha_i}y_j^0(t)=\Phi_{j+1}(t)-I_{0+}^{\alpha_0}\sum_{i=1}^{m}a_i(t)^cD_{0+}^{\alpha_i}\Phi_{j+1}(t).$$

Here

$${}^cD_{0+}^{\alpha_i}\Phi_{j+1}(t)=D_{0+}^{\alpha_i}\left[\Phi_{j+1}(t)-\sum_{k=0}^{n_i-1}D^k\Phi_{j+1}(0)\Phi_{k+1}(t)\right],\ i=1,\cdots,m.$$

Note that for $k=0,1,\cdots,n_i-1$

$$D^k\Phi_{j+1}(0)=\begin{cases}1, & k=j,\\ 0, & k\neq j.\end{cases}$$

Since $n_0>n_1\geq n_i$ then for $j=0,1,\cdots,n_1-1$, we have

$${}^cD_{0+}^{\alpha_i}\Phi_{j+1}(t)=\begin{cases}D_{0+}^{\alpha_i}\Phi_{j+1}(t), & h_j\leq i\leq m,\\ 0, & 1\leq i<h_j.\end{cases}\tag{13}$$

If $j=n_1,\cdots,n_0-1$, then $k<j$ and thus we have

$${}^cD_{0+}^{\alpha_i}\Phi_{j+1}(t)=D_{0+}^{\alpha_i}\Phi_{j+1}(t),\ \ \mathrm{i}=1,\cdots,m.\tag{14}$$

Therefore the first approximation of $y_j(t)$ is given by

$$y_j^1(t)=\Phi_{j+1}(t)-I_{0+}^{\alpha_0}\sum_{i=h_j}^{m}a_i(t)D_{0+}^{\alpha_i}\Phi_{j+1}(t)=\Phi_{j+1}(t)-I_{0+}^{\alpha_0}\sum_{i=h_j}^{m}a_i(t)\Phi_{j+1-\alpha_i}(t),$$
$$j=0,1,\cdots,n_1-1,$$

$$y_j^1(t)=\Phi_{j+1}(t)-I_{0+}^{\alpha_0}\sum_{i=1}^{m}a_i(t)D_{0+}^{\alpha_i}\Phi_{j+1}(t)=\Phi_{j+1}(t)-I_{0+}^{\alpha_0}\sum_{i=1}^{m}a_i(t)\Phi_{j+1-\alpha_i}(t),$$
$$j=n_1,\cdots,n_0-1.$$

and we have $y_j^1(t)\in C^{\alpha_0,n_0-1}[0,T],\ j=0,1,\cdots,n_0-1$.

For $j=0,1,\cdots,n_1-1$ the second approximation is given as follows:

$$y_j^2(t)=\Phi_{j+1}(t)-I_{0+}^{\alpha_0}\sum_{i=1}^{m}a_i(t)^cD_{0+}^{\alpha_i}y_j^1(t)=$$

$$=\Phi_{j+1}(t)-I_{0+}^{\alpha_0}\sum_{i=1}^{m}a_i(t)^cD_{0+}^{\alpha_i}\left[\Phi_{j+1}(t)-I_{0+}^{\alpha_0}\sum_{i=h_j}^{m}a_i(t)\Phi_{j+1-\alpha_i}(t)\right]$$

$$=\Phi_{j+1}(t)-I_{0+}^{\alpha_0}\sum_{i=1}^{m}a_i(t)^cD_{0+}^{\alpha_i}\Phi_{j+1}(t)+I_{0+}^{\alpha_0}\sum_{i=1}^{m}a_i(t)^cD_{0+}^{\alpha_i}I_{0+}^{\alpha_0}\sum_{i=h_j}^{m}a_i(t)\Phi_{j+1-\alpha_i}(t).\tag{15}$$

By (13) and the fact that ${}^cD_{0+}^{\alpha_i} I_{0+}^{\alpha_0} \sum_{i=h_j}^{m} a_i(t)\Phi_{j+1-\alpha_i}(t) = I_{0+}^{\alpha_0-\alpha_i} \sum_{i=h_j}^{m} a_i(t)\Phi_{j+1-\alpha_i}(t)$, we can rewrite (15) as follows:

$$y_j^2(t) = \Phi_{j+1}(t) - I_{0+}^{\alpha_0} \sum_{i=h_j}^{m} a_i(t) D_{0+}^{\alpha_i} \Phi_{j+1}(t) + I_{0+}^{\alpha_0} \sum_{i=1}^{m} a_i(t) I_{0+}^{\alpha_0-\alpha_i} \sum_{i=h_j}^{m} a_i(t)\Phi_{j+1-\alpha_i}(t)$$

$$= \Phi_{j+1}(t) - I_{0+}^{\alpha_0} \sum_{i=h_j}^{m} a_i(t)\Phi_{j+1-\alpha_i}(t) + I_{0+}^{\alpha_0} \sum_{i=1}^{m} a_i(t) I_{0+}^{\alpha_0-\alpha_i} \sum_{i=h_j}^{m} a_i(t)\Phi_{j+1-\alpha_i}(t)$$

$$= \Phi_{j+1}(t) + \sum_{k=0}^{1} (-1)^{k+1} I_{0+}^{\alpha_0} \left[ \sum_{i=1}^{m} a_i(t) I_{0+}^{\alpha_0-\alpha_i} \right]^k \sum_{i=h_j}^{m} a_i(t)\Phi_{j+1-\alpha_i}(t).$$

If $j = n_1, n_1+1, \cdots, n_0-1$, then the second approximation is given as follows:

$$y_j^2(t) = \Phi_{j+1}(t) - I_{0+}^{\alpha_0} \sum_{i=1}^{m} a_i(t)\, {}^cD_{0+}^{\alpha_i} y_j^1(t) =$$

$$= \Phi_{j+1}(t) - I_{0+}^{\alpha_0} \sum_{i=1}^{m} a_i(t)\, {}^cD_{0+}^{\alpha_i} \left[ \Phi_{j+1}(t) - I_{0+}^{\alpha_0} \sum_{i=1}^{m} a_i(t)\Phi_{j+1-\alpha_i}(t) \right]$$

$$= \Phi_{j+1}(t) + \sum_{k=0}^{1} (-1)^{k+1} I_{0+}^{\alpha_0} \left[ \sum_{i=1}^{m} a_i(t) I_{0+}^{\alpha_0-\alpha_i} \right]^k \sum_{i=1}^{m} a_i(t)\Phi_{j+1-\alpha_i}(t)$$

Thus we have the representation of the second approximation of $y_j(t)$:

$$y_j^2(t) = \Phi_{j+1}(t) + \sum_{k=0}^{1} (-1)^{k+1} I_{0+}^{\alpha_0} \left[ \sum_{i=1}^{m} a_i(t) I_{0+}^{\alpha_0-\alpha_i} \right]^k \sum_{i=h_j}^{m} a_i(t)\Phi_{j+1-\alpha_i}(t),$$
$$j = 0, 1, \cdots, n_1 - 1,$$

$$y_j^2(t) = \Phi_{j+1}(t) + \sum_{k=0}^{1} (-1)^{k+1} I_{0+}^{\alpha_0} \left[ \sum_{i=1}^{m} a_i(t) I_{0+}^{\alpha_0-\alpha_i} \right]^k \sum_{i=1}^{m} a_i(t)\Phi_{j+1-\alpha_i}(t),$$
$$j = n_1, n_1+1, \cdots, n_0 - 1.$$

For any $n = 3, 4, \cdots$, we can find $y_j^n(t)$ by induction and by letting $n \to \infty$ we get the representation (11) and (12) of the canonical fundamental system . It is obvious from lemma 2 that $y_j(t) \in C^{\alpha_0, n_0-1}[0,T],\ j = 0, 1, \cdots, n_0-1$ .(QED)

**Corollary 2.** *Under the assumption of theorem 2, a solution* $y(t) \in C^{\alpha_0, n_0-1}[0,T]$ *to the initial value problem (1) and (2) uniquely exists and represents by*

$$y(t) = \sum_{j=0}^{n_0-1} b_j y_j(t)$$

*Here* $y_j(t) \in C^{\alpha_0, n_0 - 1}[0,T],\ j = 0, 1, \cdots, n_0 - 1$ *is the canonical fundamental system of (1) given by (11) and (12).*

## 4. The Canonical Fundamental System of Solutions to the Equations of the Type II

**Theorem 3.** *Let* $0 < \gamma = n_0 - \alpha_0 < 1$ *and assume that* $n_0 = n_1$ *and* $a_i \in C_\gamma^1[0,T]$ , $D_{0+}^{n_0 - \alpha_0} a_i \in C[0,T], i = 1, \cdots, m$ . *Assume that the equation (1)is the type II. That is, assume that* $n_0 \geq 2$ and there exists a $j_0 \in \{0, 1, \cdots, n_0 - 2\}$ such that $H_{j_0} = \phi$ and $H_{j_0+1} \neq \phi$ . *Then there exists the unique canonical fundamental system* $y_j(t) \in C_\gamma^{n_0}[0,T]$ , $j = 0, 1, \cdots, n_0 - 1$ *of solutions to (1) and it is written by*

$$y_j(t) = \Phi_{j+1}(t), \qquad j = 0, 1, \cdots, j_0, \tag{16}$$

$$y_j(t) = \Phi_{j+1}(t) + \sum_{k=0}^{\infty} (-1)^{k+1} I_{0+}^{\alpha_0} \left[ \sum_{i=1}^{m} a_i(t) I_{0+}^{\alpha_0 - \alpha_i} \right]^k \sum_{i=h_j}^{m} a_i(t) \Phi_{j+1-\alpha_i}(t), \tag{17}$$
$$j = j_0 + 1, \cdots, n_0 - 1.$$

**Proof**. From lemma 2-(ii) the existence and uniqueness of canonical fundamental systems is evident. For $j \in \{0, 1, \cdots, n_0 - 1\}$ , we find the $j$ -th element $y_j(t)$ of the of canonical fundamental systems as the limit in $C_\gamma^{n_0}[0,T]$ of the approximation sequence :

$$y_j^0(t) = \Phi_{j+1}(t) \tag{18}$$

$$y_j^l(t) = y_j^0(t) - I_{0+}^{\alpha_0} \sum_{i=1}^{m} a_i(t) {}^c D_{0+}^{\alpha_i} y_j^{l-1}(t),\ l = 1, 2, \cdots. \tag{19}$$

The first approximation of $y_j(t)$ for $j = 0, 1, \cdots, n_0 - 1$ can be written by

$$y_j^1(t) = \Phi_{j+1}(t) - I_{0+}^{\alpha_0} \sum_{i=1}^{m} a_i(t) {}^c D_{0+}^{\alpha_i} y_j^0(t) = \Phi_{j+1}(t) - I_{0+}^{\alpha_0} \sum_{i=1}^{m} a_i(t) {}^c D_{0+}^{\alpha_i} \Phi_{j+1}(t). \tag{20}$$

We need to calculate ${}^c D_{0+}^{\alpha_i} \Phi_{j+1}(t)$ .

$${}^c D_{0+}^{\alpha_i} \Phi_{j+1}(t) = D_{0+}^{\alpha_i} \left[ \Phi_{j+1}(t) - \sum_{k=0}^{n_i - 1} D^k \Phi_{j+1}(0) \Phi_{k+1}(t) \right],\ i = 1, \cdots, m. \tag{21}$$

Here note that

$$D^k \Phi_{j+1}(0) = \begin{cases} 1, & k = j, \\ 0, & k \neq j, \end{cases}$$

and thus

$$D^k \Phi_{j+1}(0) \Phi_{k+1}(t) = \begin{cases} \Phi_{j+1}(t), & k = j, \\ 0, & k \neq j. \end{cases}$$

Therefore in order to calculate ${}^cD_{0+}^{\alpha_i}\Phi_{j+1}(t)$ by (21), for every fixed $j=0,1,\cdots,n_0-1$ and $i=1,\cdots,m$, we must consider the relation between $j$ and $k=0,\ 1,\cdots,n_i-1$. From the assumption there exists a $j_0\in\{0,\ 1,\cdots,n_0-2\}$ such that

$$H_{j_0}=\{i:0\le\alpha_i\le j_0,\ i=1,\cdots,m\}=\phi \text{ and } H_{j_0+1}=\{i:0\le\alpha_i\le j_0+1,\ i=1,\cdots,m\}\ne\phi.$$

Thus $\alpha_i>j_0$ for every $i=1,2,\cdots,m$ and from $n_i-1<\alpha_i\le n_i$ we have $j_0\le n_i-1$.

Therefore if $j=0,1,\cdots,j_0$, for all $i\in\{1,\cdots,m\}$ there exists a $k\in\{0,1,\cdots,n_i-1\}$ such that $k=j$ and thus

$${}^cD_{0+}^{\alpha_i}\Phi_{j+1}(t)=D_{0+}^{\alpha_i}\left[\Phi_{j+1}(t)-\sum_{k=0}^{n_i-1}D^k\Phi_{j+1}(0)\Phi_{k+1}(t)\right]=0,\ i=1,2,\cdots,m. \tag{22}$$

Next, for $j=j_0+1,\cdots,n_0-1$, we calculate ${}^cD_{0+}^{\alpha_i}\Phi_{j+1}(t)$. Let $j=j_0+1$. Note that

$$j_0<\alpha_m<\alpha_{m-1}<\cdots<\alpha_1\le n_1$$

and

$$h_{j_0+1}=\min H_{j_0+1}=\min\{i:0\le\alpha_i\le j_0+1,\ i=1,\cdots,m\}.$$

Then for all $i=h_{j_0+1},h_{j_0+1}+1,\cdots,m$, we have $j_0<\alpha_i\le j_0+1$. On the other hand, since $n_i-1<\alpha_i\le n_i$, we have $j_0=n_i-1\,(\forall i=h_{j_0+1},h_{j_0+1}+1,\cdots,m)$ and thus for any $k=0,1,\cdots,n_i-1$, we have $k<j_0+1=j$. For any $i=h_{j_0+1},h_{j_0+1}+1,\cdots,m$, we have

$$\begin{aligned}{}^cD_{0+}^{\alpha_i}\Phi_{j+1}(t)&={}^cD_{0+}^{\alpha_i}\Phi_{j_0+2}(t)=D_{0+}^{\alpha_i}\left[\Phi_{j_0+2}(t)-\sum_{k=0}^{n_i-1}D^k\Phi_{j_0+2}(0)\Phi_{k+1}(t)\right]=D_{0+}^{\alpha_i}\Phi_{j_0+2}(t)\\&=\Phi_{j_0+2-\alpha_i}(t)=\Phi_{j+1-\alpha_i}(t).\end{aligned} \tag{23}$$

Note that $j_0+1<\alpha_{h_{j_0+1}-1}<\cdots<\alpha_1$ and $n_i-1<\alpha_i\le n_i$. Thus for $i\in\{1,2,\cdots,h_{j_0+1}-1\}$, we have $j=j_0+1\le n_i-1$ and

$$\begin{aligned}{}^cD_{0+}^{\alpha_i}\Phi_{j+1}(t)&={}^cD_{0+}^{\alpha_i}\Phi_{j_0+2}(t)=D_{0+}^{\alpha_i}\left[\Phi_{j_0+2}(t)-\sum_{k=0}^{n_i-1}D^k\Phi_{j_0+2}(0)\Phi_{k+1}(t)\right]=\\&=D_{0+}^{\alpha_i}[\Phi_{j_0+2}(t)-\Phi_{j_0+2}(t)]=0.\end{aligned} \tag{24}$$

Putting (23) and (24) together, when $j=j_0+1$, we have

$${}^cD_{0+}^{\alpha_i}\Phi_{j+1}(t)=\begin{cases}\Phi_{j+1-\alpha_i}(t), & i=h_j,h_j+1,\cdots,m,\\ 0, & i=1,2,\cdots,h_j-1.\end{cases} \tag{25}$$

Next more generally, let $j=j_0+l$, where $l=\{1,\cdots,n_0-j_0-1\}$. From the definition of $h_j$, for any $i\in\{h_j,h_j+1,\cdots,m\}$, we have $\alpha_i\le j$. On the other hand, from $n_i-1<\alpha_i\le n_i$, we have $n_i\le j$ and thus for any $k=0,1,\cdots,n_i-1$, we have $k<j$ and

$${}^cD_{0+}^{\alpha_i}\Phi_{j+1}(t)=D_{0+}^{\alpha_i}\left[\Phi_{j+1}(t)-\sum_{k=0}^{n_i-1}D^k\Phi_{j+1}(0)\Phi_{k+1}(t)\right]$$

$$= D_{0+}^{\alpha_i}\Phi_{j+1}(t) = \Phi_{j+1-\alpha_i}(t).$$

If $i \in \{1, 2, \cdots, h_j - 1\}$ , then $j < \alpha_i$ and thus $j \le n_i - 1$ . Thus we have ${}^cD_{0+}^{\alpha_i}\Phi_{j+1}(t) = 0.$ Putting these together, for $j = j_0 + 1, \cdots, n_0 - 1$, we get

$$ {}^cD_{0+}^{\alpha_i}\Phi_{j+1}(t) = \begin{cases} \Phi_{j+1-\alpha_i}(t), & i = h_j, \cdots, m, \\ 0, & i = 1, 2, \cdots, h_j - 1. \end{cases} \tag{26}$$

From (20), (22) and (26), the first approximation of $y_j(t)$ is written by

$$\left. \begin{aligned} & y_j^1(t) = \Phi_{j+1}(t), && j = 0, 1, \cdots, j_0, \\ & y_j^1(t) = \Phi_{j+1}(t) - I_{0+}^{\alpha_0} \sum_{i=h_j}^m a_i(t)\Phi_{j+1-\alpha_i}(t), && j = j_0 + 1, \cdots, n_0 - 1 \end{aligned} \right\} \tag{27}$$

and $y_j^1(t) \in C_\gamma^{n_0}[0,T], \ j = j_0 + 1, \cdots, n_0 - 1$ .

Next we find the second approximation of $y_j(t)$. From (18) and (19) we can write

$$y_j^2(t) = \Phi_{j+1}(t) - I_{0+}^{\alpha_0} \sum_{i=1}^m a_i(t)\, {}^cD_{0+}^{\alpha_i} y_j^1(t), \quad j = 0, 1, \cdots, n_0 - 1$$

and by (27) we have

$$y_j^2(t) = \Phi_{j+1}(t) - I_{0+}^{\alpha_0} \sum_{i=1}^m a_i(t)\, {}^cD_{0+}^{\alpha_i}\Phi_{j+1}(t), \quad j = 0, 1, \cdots, j_0.$$

From (22), ${}^cD_{0+}^{\alpha_i}\Phi_{j+1}(t) = 0, \ i = 1, \cdots, m$ and thus

$$y_j^2(t) = \Phi_{j+1}(t), \qquad j = 0, 1, \cdots, j_0.$$

Let $j \in \{j_0 + 1, \cdots, n_0 - 1\}$ . From (27), we can write

$$y_j^2(t) = \Phi_{j+1}(t) - I_{0+}^{\alpha_0} \sum_{i=1}^m a_i(t)\, {}^cD_{0+}^{\alpha_i} \left[ \Phi_{j+1}(t) - I_{0+}^{\alpha_0} \sum_{i=h_j}^m a_i(t)\Phi_{j+1-\alpha_i}(t) \right]$$

$$= \Phi_{j+1}(t) - I_{0+}^{\alpha_0} \sum_{i=1}^m a_i(t)\, {}^cD_{0+}^{\alpha_i}\Phi_{j+1}(t) + I_{0+}^{\alpha_0} \sum_{i=1}^m a_i(t)\, {}^cD_{0+}^{\alpha_i} I_{0+}^{\alpha_0} \sum_{i=h_j}^m a_i(t)\Phi_{j+1-\alpha_i}(t).$$

By (26), we have

$${}^cD_{0+}^{\alpha_i}\Phi_{j+1}(t) = \Phi_{j+1-\alpha_i}(t), \quad j = j_0 + 1, \cdots, n_0 - 1, \ i = h_j, \cdots, m.$$

Note that ${}^cD_{0+}^{\alpha_i} I_{0+}^{\alpha_0} \sum_{i=h_j}^m a_i(t)\Phi_{j+1-\alpha_i}(t) = I_{0+}^{\alpha_0-\alpha_i} \sum_{i=h_j}^m a_i(t)\Phi_{j+1-\alpha_i}(t)$. Then we get

$$y_j^2 = \Phi_{j+1}(t) - I_{0+}^{\alpha_0} \sum_{i=h_j}^m a_i(t)\Phi_{j+1-\alpha_i}(t) + I_{0+}^{\alpha_0} \sum_{i=1}^m a_i(t) I_{0+}^{\alpha_0-\alpha_i} \sum_{i=h_j}^m a_i(t)\Phi_{j+1-\alpha_i}(t)$$

$$= \Phi_{j+1}(t) + \sum_{k=0}^{1} (-1)^k I_{0+}^{\alpha_0} \left[ \sum_{i=1}^m a_i(t) I_{0+}^{\alpha_0-\alpha_i} \right]^k \sum_{i=h_j}^m a_i(t)\Phi_{j+1-\alpha_i}(t).$$

That is, the second approximation of $y_j(t)$ is given as follows:

$$y_j^2(t) = \Phi_{j+1}(t), \qquad j = 0, 1, \cdots, j_0;$$

$$y_j^2 = \Phi_{j+1}(t) + \sum_{k=0}^{1} (-1)^k I_{0+}^{\alpha_0} \left[ \sum_{i=1}^{m} a_i(t) I_{0+}^{\alpha_0 - \alpha_i} \right]^k \sum_{i=h_j}^{m} a_i(t) \Phi_{j+1-\alpha_i}(t),$$

$$j = j_0 + 1, \cdots, n_0 - 1.$$

By induction we get the $n$-th approximation of $y_j(t)$:

$$y_j^n(t) = \Phi_{j+1}(t), \qquad j = 0, 1, \cdots, j_0;$$

$$y_j^n = \Phi_{j+1}(t) + \sum_{k=0}^{n-1} (-1)^{k+1} I_{0+}^{\alpha_0} \left[ \sum_{i=1}^{m} a_i(t) I_{0+}^{\alpha_0 - \alpha_i} \right]^k \sum_{i=h_j}^{m} a_i(t) \Phi_{j+1-\alpha_i}(t),$$

$$j = j_0 + 1, \cdots, n_0 - 1.$$

Thus the canonical fundamental system of solutions to (1) is given by (16) and (17), and we have $y_j(t) \in C_\gamma^{n_0}[0, T],\ j = 0, 1, \cdots, n_0 - 1$ . (QED)

**Corollary 3.** *Under the assumption of theorem 3, a solution* $y(t) \in C_\gamma^{n_0}[0, T]$ *to the initial value problem (1) and (2) uniquely exists and represents by* $y(t) = \sum_{j=0}^{n_0 - 1} b_j y_j(t)$. Here $y_j(t) \in C_\gamma^{n_0}[0, T],\ j = 0, 1, \cdots, n_0 - 1$ *is the canonical fundamental system of (1) given by* (16) and (17).

**Theorem 4.** *Assume that* $n_0 > n_1$ *and* $a_i \in C[0, T], i = 1, \cdots, m$. *Assume that the equation (1)is the type II. That is, assume that* $n_0 \geq 2$ and there exists a $j_0 \in \{0, 1, \cdots, n_0 - 2\}$ such that $H_{j_0} = \phi$ and $H_{j_0 + 1} \neq \phi$. *Then there exists the unique canonical fundamental system* $y_j(t) \in C^{\alpha_0, n_0 - 1}[0, T],\ j = 0, 1, \cdots, n_0 - 1$ *of solutions to (1) and it is written by*

$$y_j(t) = \Phi_{j+1}(t), \qquad j = 0, 1, \cdots, j_0, \tag{28}$$

$$y_j(t) = \Phi_{j+1}(t) + \sum_{k=0}^{\infty} (-1)^{k+1} I_{0+}^{\alpha_0} \left[ \sum_{i=1}^{m} a_i(t) I_{0+}^{\alpha_0 - \alpha_i} \right]^k \sum_{i=h_j}^{m} a_i(t) \Phi_{j+1-\alpha_i}(t), \tag{29}$$

$$j = j_0 + 1, \cdots, n_1 - 1,$$

$$y_j(t) = \Phi_{j+1}(t) + \sum_{k=0}^{\infty} (-1)^{k+1} I_{0+}^{\alpha_0} \left[ \sum_{i=1}^{m} a_i(t) I_{0+}^{\alpha_0 - \alpha_i} \right]^k \sum_{i=1}^{m} a_i(t) \Phi_{j+1-\alpha_i}(t), \tag{30}$$

$$j = n_1, n_1 + 1, \cdots, n_0 - 1.$$

**Proof.** From lemma 2(i) the existence and uniqueness of canonical fundamental systems is evident. For $j \in \{0, 1, \cdots, n_0 - 1\}$, we find the $j$-th element $y_j(t)$ of the of canonical fundamental systems as the limit of the approximation sequence (18) and (19). Let $y_j^0(t) = \Phi_{j+1}(t)$, then by (19), the first approximation of $y_j(t)$ is given as follows:

$$y_j^1(t) = \Phi_{j+1}(t) - I_{0+}^{\alpha_0} \sum_{i=1}^{m} a_i(t) {}^c D_{0+}^{\alpha_i} y_j^{0`}(t)$$

$$= \Phi_{j+1}(t) - I_{0+}^{\alpha_0} \sum_{i=1}^{m} a_i(t) {}^c D_{0+}^{\alpha_i} \Phi_{j+1}(t) \ , \ j = 0, 1, \cdots, n_0 - 1 .$$

Here calculate ${}^c D_{0+}^{\alpha_i} \Phi_{j+1}(t)$. From the definition 2,

$${}^c D_{0+}^{\alpha_i} \Phi_{j+1}(t) = D_{0+}^{\alpha_i} \left[ \Phi_{j+1}(t) - \sum_{k=0}^{n_i - 1} D^k \Phi_{j+1}(0) \Phi_{k+1}(t) \right], \ i = 1, \cdots, m \tag{31}$$

We know

$$D^k \Phi_{j+1}(0) \Phi_{k+1}(t) = \begin{cases} \Phi_{j+1}(t), & k = j, \\ 0, & k \neq j. \end{cases}$$

From the assumption $H_{j_0} = \phi$ and $H_{j_0+1} \neq \phi$, we have $\alpha_i > j_0$ for all $i = 1,2, \cdots, m$. On the other hand, $n_i - 1 < \alpha_i \le n_i$ and thus $j_0 \le n_i - 1$. Therefore for any $j = 0, 1, \cdots, j_0$ there exists a $k \in \{0, \cdots n_i - 1\}$ such that $j = k$. Thus

$${}^c D_{0+}^{\alpha_i} \Phi_{j+1}(t) = D_{0+}^{\alpha_i} [\Phi_{j+1}(t) - \Phi_{j+1}(t)] = 0 \ , \ j = 0, 1, \cdots, j_0$$

and thus

$$y_j^1(t) = \Phi_{j+1}(t) , \quad j = 0, 1, \cdots, j_0 . \tag{32}$$

For $j = j_0 + 1, \cdots, n_1 - 1$, we get

$${}^c D_{0+}^{\alpha_i} \Phi_{j+1}(t) = \begin{cases} \Phi_{j+1-\alpha_i}(t), & i = h_j, \cdots, m , \\ 0 \quad , & i = 1,2, \cdots, h_j - 1. \end{cases}$$

as the same consideration with (26) and thus we have.

$$y_j^1(t) = \Phi_{j+1}(t) - I_{0+}^{\alpha_0} \sum_{i=h_j}^{m} a_i(t) \Phi_{j+1-\alpha_i}(t), \quad j = j_0 + 1, \cdots, n_1 - 1 . \tag{33}$$

For any $j = n_1, n_1 + 1, \cdots, n_0 - 1$ , since $n_1 \ge n_2 \ge \cdots \ge n_m$ , we always have $k < j$ when $k = 0, 1, \cdots, n_i - 1 (i = 1, \cdots, m)$. Thus we always have $D^k \Phi_{j+1}(0) \Phi_{k+1}(t) = 0$ and therefore

$${}^c D_{0+}^{\alpha_i} \Phi_{j+1}(t) = D_{0+}^{\alpha_i} \Phi_{j+1}(t) = \Phi_{j+1-\alpha_i}(t), \ \ j = n_1, n_1 + 1, \cdots, n_0 - 1 .$$

So we get

$$y_j^1(t) = \Phi_{j+1}(t) - I_{0+}^{\alpha_0} \sum_{i=1}^{m} a_i(t) \Phi_{j+1-\alpha_i}(t) \tag{34}$$

Put (32),(33) and (34) together, we have

$$y_j^1(t) = \Phi_{j+1}(t) , \qquad j = 0, 1, \cdots, j_0 ,$$

$$y_j^1(t) = \Phi_{j+1}(t) - I_{0+}^{\alpha_0} \sum_{i=h_j}^{m} a_i(t)\Phi_{j+1-\alpha_i}(t), \quad j = j_0+1, \cdots, n_1-1,$$

$$y_j^1(t) = \Phi_{j+1}(t) - I_{0+}^{\alpha_0} \sum_{i=1}^{m} a_i(t)\Phi_{j+1-\alpha_i}(t), \quad j = n_1, n_1+1, \cdots, n_0-1,$$

In similar way as the above we get the second approximation of $y_j(t)$:

$$y_j^2(t) = \Phi_{j+1}(t), \qquad j = 0, 1, \cdots, j_0,$$

$$y_j^2(t) = \Phi_{j+1}(t) + \sum_{k=0}^{1} (-1)^k I_{0+}^{\alpha_0} \left[ \sum_{i=1}^{m} a_i(t) I_{0+}^{\alpha_0-\alpha_i} \right]^k \sum_{i=h_j}^{m} a_i(t)\Phi_{j+1-\alpha_i}(t), \quad j = j_0+1, \cdots, n_1-1,$$

$$y_j^2(t) = \Phi_{j+1}(t) + \sum_{k=0}^{1} (-1)^k I_{0+}^{\alpha_0} \left[ \sum_{i=1}^{m} a_i(t) I_{0+}^{\alpha_0-\alpha_i} \right]^k \sum_{i=1}^{m} a_i(t)\Phi_{j+1-\alpha_i}(t), \quad j = n_1, n_1+1, \cdots, n_0-1.$$

By induction we can get the $n$-th approximation of $y_j(t)$ and then taking $n \to \infty$ we get the canonical fundamental system

$$y_j(t) = \Phi_{j+1}(t), \qquad j = 0, 1, \cdots, j_0,$$

$$y_j(t) = \Phi_{j+1}(t) + \sum_{k=0}^{\infty} (-1)^{k+1} I_{0+}^{\alpha_0} [\sum_{i=1}^{m} a_i(t) I_{0+}^{\alpha_0-\alpha_i}]^k \sum_{i=h_j}^{m} a_i(t)\Phi_{j+1-\alpha_i}(t), \quad j = j_0+1, \cdots, n_1-1,$$

$$y_j(t) = \Phi_{j+1}(t) + \sum_{k=0}^{\infty} (-1)^{k+1} I_{0+}^{\alpha_0} [\sum_{i=1}^{m} a_i(t) I_{0+}^{\alpha_0-\alpha_i}]^k \sum_{i=1}^{m} a_i(t)\Phi_{j+1-\alpha_i}(t), \quad j = n_1, n_1+1, \cdots, n_0-1$$

and $y_j(t) \in C^{\alpha_0, n_0-1}[0,T],\ j = 0, 1, \cdots, n_0-1$. (QED)

**Corollary 4.** *Under the assumption of theorem 4, a solution* $y(t) \in C^{\alpha_0, n_0-1}[0,T]$ *to the initial value problem (1) and (2) uniquely exists and represents by* $y(t) = \sum_{j=0}^{n_0-1} b_j y_j(t)$.

Here $y_j(t) \in C^{\alpha_0, n_0-1}[0,T],\ j = 0, 1, \cdots, n_0-1$ *is the canonical fundamental system of (1) given by* (28) (29) and (30).

## 5. The Canonical Fundamental System of Solutions to the Equations of the Type III

**Theorem 5.** *Let* $0 < \gamma = n_0 - \alpha_0 < 1$ *and assume that* $a_i \in C_\gamma^1[0,T], i = 1, \cdots, m$. *Assume that the equation (1)is the type III. That is, assume that* $H_{n_0-1} = \phi$. *Then (1) has the unique canonical fundamental system* $\{y_j(t): j = 0, 1, \cdots, n_0-1\}$ *in* $C_\gamma^{n_0}[0,T]$ *of solutions and it is represented as follows:*

$$y_j(t) = \Phi_{j+1}(t), \ \ j = 0, 1, \cdots, n_0 - 1. \tag{35}$$

**Proof.** From lemma 2(ii) the existence and uniqueness of canonical fundamental system in $C_\gamma^{n_0}[0,T]$ is evident. Let $y_j^0(t) = \Phi_{j+1}(t)$ , $j = 0, 1, \cdots, n_0 - 1$ and calculate the first approximation of $y_j(t)$, the we have

$$y_j^1 = \Phi_{j+1}(t) - I_{0+}^{\alpha_0} \sum_{i=1}^{m} a_i(t) {}^c D_{0+}^{\alpha_i} \Phi_{j+1}(t).$$

From the assumption of $H_{n_0-1} = \{i : 0 \le \alpha_i \le n_0 - 1, \ i = 1, \cdots, m\} = \phi$ , it is evident that $n_0 - 1 < \alpha_i \le n_0$ for all $i = 1, \cdots, m$ and thus

$${}^c D_{0+}^{\alpha_i} \Phi_{j+1}(t) = D_{0+}^{\alpha_i} [\Phi_{j+1}(t) - \sum_{k=0}^{n_0-1} D^k \Phi_{j+1}(0) \Phi_{j+1}(t)] = 0, \ \ j = 0, 1, \cdots, n_0 - 1$$

and $y_j^1(t) = \Phi_{j+1}(t)$, $j = 0, 1, \cdots, n_0 - 1$. Similarly, we easily know

$$y_j^n(t) = \Phi_{j+1}(t), \ \ j = 0, 1, \cdots, n_0 - 1$$

and thus we get

$$y_j(t) = \Phi_{j+1}(t) \in C^\infty[0,T], \ \ j = 0, 1, \cdots, n_0 - 1.$$

(QED)

**Corollary 5.** *Under the assumption of theorem 5, a solution* $y(t) \in C_\gamma^{n_0}[0,T]$ *to the initial value problem (1) and (2) uniquely exists and represents by* $y(t) = \sum_{j=0}^{n_0-1} b_j y_j(t)$. Here $y_j(t) \in C^\infty[0,T], \ j = 0, 1, \cdots, n_0 - 1$ *is the canonical fundamental system of (1) given by* (35).

***Remark 4.*** The theorem 5 and corollary 5 show that the equations of type III have the solution independent of the coefficients $a_i(t)$ and orders $\alpha_i$. It reflects the essence of the equations of type III.

## 6. Conclusions

First, the canonical fundamental system of solutions to linear homogeneous differential equation (1) with Caputo fractional derivatives and continuous variable coefficients has different representations according to the distributions of the lowest order $\alpha_m$. That is, it has different representations depending on whether $\alpha_m = 0$ (case I) or $0 < \alpha_m \le n_0 - 1$ (case II) or $n_0 - 1 \le \alpha_m \le n_0$ (case III).

Second, in the case I and case II, the canonical fundamental system of solutions to linear homogeneous differential equation (1) with Caputo fractional derivatives and continuous variable coefficients has different representations according to the distance between the two highest orders $\alpha_0$ and $\alpha_1$.

Third, 5 *patterns of distributions of fractional orders* in the equation (1) which *determine the representations* of the solutions are as follows:

1) $\alpha_m = 0$ and $n_0 = n_1$ (theorem 1 and corollary 1) ;
2) $\alpha_m = 0$ and $n_0 > n_1$ (theorem 2 and corollary 2);
3) $0 < \alpha_m \le n_0 - 1$ and $n_0 = n_1$ (theorem 3 and corollary 3) ;
4) $0 < \alpha_m \le n_0 - 1$ and $n_0 > n_1$ (theorem 4 and corollary 4) ;
5) $n_0 - 1 \le \alpha_m \le n_0$ (theorem 5 and corollary 5).